\documentclass[letterpaper, 10 pt, conference]{ieeeconf}  
\IEEEoverridecommandlockouts                              

\let\labelindent\relax
\usepackage{enumitem}
\usepackage{booktabs}
\usepackage{graphicx}
\usepackage{amssymb,mathtools,amsfonts}
\usepackage{gensymb}
 
\usepackage{amsthm}
\usepackage{caption, subcaption}
\usepackage{hyperref}
\usepackage{cleveref}
\usepackage{epsfig}
\usepackage{multirow}
\usepackage{color}
\usepackage{psfrag}
\usepackage{empheq}
\usepackage{bm} 
\usepackage{url}
\usepackage{float}
\usepackage{comment}
\usepackage{alphalph}
\usepackage[short]{optidef}
\usepackage[dvipsnames]{xcolor}
\usepackage{tikz}
\usepackage{tikz-cd}
\usepackage{suffix}

\usetikzlibrary{positioning}
\usetikzlibrary{arrows}
\usetikzlibrary{babel}
\usetikzlibrary{automata}
\usetikzlibrary{shapes, snakes}
\usetikzlibrary{positioning, fit, calc}
\usetikzlibrary{decorations.pathreplacing}
\usetikzlibrary{patterns}
\usetikzlibrary{patterns.meta}
\usetikzlibrary{arrows.meta}
\usepackage{pgfplots}
\pgfplotsset{compat=newest}
\usepackage{pgfplotstable}
\usepackage{pgf}
\usepackage{cite}

\usepackage[export]{adjustbox}

\newcommand{\R}{{\mathbb{R}}}

\newcommand{\nstate}{n_\mathrm{x}}
\newcommand{\nalg}{n_\mathrm{z}}
\newcommand{\ncontrol}{n_\mathrm{u}}
\newcommand{\ns}{n_\mathrm{s}}

\newcommand{\paren}[1]{\left( #1 \right)}
\newcommand{\brac}[1]{\left[ #1 \right]}
\newcommand{\cbrac}[1]{\left\{ #1 \right\}}
\WithSuffix\newcommand\paren*[1]{( #1 )}
\WithSuffix\newcommand\brac*[1]{[ #1 ]}
\WithSuffix\newcommand\cbrac*[1]{\{ #1 \}}

\newcommand{\vo}{\vec{o}\@ifnextchar{^}{\,}{}}

\newcommand{\Set}[2]{\left\{\, #1 \mid #2\,\right\}}
\WithSuffix\newcommand\Set*[2]{\{\, #1 \mid #2\,\}}

\newcommand{\transp}[1]{{#1}^\top}
\newcommand{\inv}[1]{{#1}^{-1}}

\newcommand{\xdot}{{\dot{x}}}
\newcommand{\proj}[2]{\mathrm{P}_{#1}\paren{#2}}
\WithSuffix\newcommand\proj*[2]{\mathrm{P}_{#1}\paren*{#2}}

\newcommand{\acados}{\texttt{acados}}
\newcommand{\blasfeo}{\texttt{BLASFEO}}
\newcommand{\fatrop}{\texttt{fatrop}}

\newcommand{\cnewton}{\zeta_{\textrm{newton}}}
\newcommand{\cfactor}{\zeta_{\textrm{factor}}}
\newcommand{\cbacksolve}{\zeta_{\textrm{backsolve}}}
\newcommand{\csens}{\zeta_{\textrm{sens}}}
\newcommand{\ctotal}{\zeta_{\textrm{total}}}

\newcommand{\cfactorhat}{\hat{\zeta}_{\textrm{factor}}}
\newcommand{\cbacksolvehat}{\hat{\zeta}_{\textrm{backsolve}}}

\newcommand{\ctotalhat}{\hat{\zeta}_{\textrm{total}}}

\begin{document}

\title{Efficient Implementation of Diagonally Implicit Runge-Kutta Methods for Model Predictive Control}
\author{
  Anton Pozharskiy$^{1}$, Jonathan Frey$^{1}$, Katrin Baumg{\"a}rtner$^{1}$, Armin Nurkanovi\'c$^{1}$, Moritz Diehl$^{1,2}$
  \thanks{This research was supported by DFG via Research Unit FOR 2401, project 424107692 and 525018088, by BMWK via 03EI4057A and 03EN3054B, and by the EU via ELO-X 953348.}
  \thanks{
    Contact: \texttt{\{first.last\}@imtek.uni-freiburg.de}
    }
  \thanks{
  $^1$Department of Microsystems Engineering (IMTEK),
  $^2$Department of Mathematics, University of Freiburg, Germany,
}
\thanks{\textbf{No} AI tools were used for the creation of this text, nor for the software implementation.}
}
\maketitle
\thispagestyle{empty} 
\begin{abstract}
  In direct shooting methods for nonlinear optimal control, one of the major computational costs is simulation of the system and calculation of sensitivities.
  These costs are particularly steep when the stiffness of the controlled system prescribes the use of implicit integrators, amongst which, implicit Runge-Kutta schemes are most popular.
  These methods require the use of Newton's method to solve a moderately sized system of nonlinear equations with a dense Jacobian, the computation of which grows cubically with the number of stages, that is, the number of abscissa.
  To mitigate this we investigate the use of the family of Diagonally Implicit Runge-Kutta (DIRK) methods whose Butcher tableaux are lower triangular.
  We describe the structure exploiting linear algebra required to make the DIRK methods competitive with fully implicit methods, including efficient sensitivity calculations.
  Finally, we evaluate the performance of DIRK methods on several example systems, and show that the efficient implementation yields Pareto optimal integration schemes.
  The fixed step size integrators are provided as part of the {\acados} package and use high-performance linear algebra provided by {\blasfeo}.
\end{abstract}
\section{Introduction}
In Nonlinear Model Predictive Control (NMPC) a feedback scheme is realized by repeatedly (and often approximately) solving the nonlinear optimal control problem:
\begin{mini}
  {\substack{x_0,\ldots, x_N\\ u_0,\ldots,u_{N-1}}}
  {\sum_{i=0}^{N-1} L(x_i, u_i) + M(x_N)}
  {}
  {}
  \addConstraint{x_0}{= \hat{x}_0}
  \addConstraint{x_{i+1}}{= F(x_i, u_i),\ }{i = 0,\ldots,N-1}
  \addConstraint{g_i(x_i, u_i)}{\le 0,\ }{i = 1,\ldots,N-1}
  \addConstraint{g_N(x_N)}{\le 0,}
\end{mini}
where $F(x_i,u_i): \R^{\nstate}\times\R^{\ncontrol} \rightarrow \R^{\nstate}$ is the numerical discretization of a continuous-time semi-explicit Differential-Algebraic Equation (DAE):
\begin{subequations}
  \label{eq:dae}
  \begin{align}
    \dot{x} &= f(x,z,u),\\
    0 &= g(x,z,u),
  \end{align}
\end{subequations}
with differential state $x\in \R^{\nstate}$, algebraic state $z\in \R^{\nalg}$, and control $u\in \R^{\ncontrol}$.
Specialized solvers for NMPC such as {\acados}~\cite{Verschueren2021} and {\fatrop}~\cite{Vanroye2023} focus on sub-problems with a linear multiple shooting structure.
A major component of building these sub-problems is the linearization of the dynamics $x_{i+1}= F(x_i, u_i)$:
\begin{equation}
  \label{eq:lin_dyn}
  x_{i+1} = A_ix_i + B_iu_i,
\end{equation}
with linear dynamics $A_i = \nabla_{x_i}F(x_i,u_i)\in \R^{\nstate\times\nstate}$, and $B_i = \nabla_{u_i}F(x_i,u_i)\in \R^{\nstate\times\ncontrol}$.


In order to formulate the linearized dynamics, numerical integration schemes which can provide sensitivities are required.
In many applications, the continuous time DAE is in fact a non-stiff Ordinary Differential Equation (ODE) and can be easily handled via explicit integration schemes such as the famous fourth order explicit scheme RK4.
However, when the dynamics are numerically stiff, such explicit schemes are inefficient and implicit schemes must be used~\cite{Hairer1989}.
In these cases it is often suggested~\cite{Rawlings2017} that practitioners use members of the Gau\ss-Legendre or Radau-IIA classes of methods, which are both Fully-Implicit Runge-Kutta (FIRK) methods.
This is because these methods are of the order proportional to twice the number of stages, and thus offer good efficiency regarding the trade-off between accuracy and evaluation cost.

The field of numerical integration using Runge-Kutta methods has, however, spent many decades developing schemes with Butcher tableaux~\cite{Butcher1976} with exploitable structure.
These include a zoo of Diagonally Implicit Runge-Kutta (DIRK) methods~\cite{Alexander1977,Cameron1983,Kennedy2016}, along with the Singly Implicit Runge-Kutta (SIRK) methods~\cite{Nilsen1997,Kennedy2025}.
These structured tableaux allow for more efficient solving of the nonlinear systems of equations by accelerating the factorization of the associated linearizations, when applying Newton's method~\cite{Kennedy2016}.
These structured methods are, however, only rarely used in NMPC contexts~\cite{Flaig2023,Schytt2024,Christensen2025}, and we are unaware of any efficient open-source implementations tailored for NMPC.
An alternative approach for structure exploitation in the Runge-Kutta equations is to detect and exploit linear dependencies in the model equations.
This is also implemented in \texttt{acados}~\cite{Frey2019} and might even be combined with a structured-tableau approach.

\begin{figure}[t]
  \centering
  \begin{tikzpicture}
  \fill[color=gray] (0,0) -- (1.7,0) -- (0,1.7) -- cycle;
  \draw[line width=2pt] (0,0) -- (2,0) -- (2,2) -- (0,2) -- cycle;
  \node (erk) at (1,-0.5) {ERK};

  \fill[color=gray] (3,0) -- (5,0) -- (5,2) -- (3,2) -- cycle;
  \draw[line width=2pt] (3,0) -- (5,0) -- (5,2) -- (3,2) -- cycle;
  \node (firk) at (4,-0.5) {FIRK};

  \fill[color=gray] (6,0) -- (7.7,0) -- (6,1.7) -- cycle;
  \draw[line width=2pt] (6,0) -- (8,0) -- (8,2) -- (6,2) -- cycle;
  \draw[line width=1] (6.3,2) -- (8,0.3);
  \node (g1) at (6.25,1.75) {\large$\gamma$};
  \node (g2) at (6.5,1.5) {\large$\gamma$};
  \node (g3) at (6.75,1.25) {\large$\gamma$};
  \node (g4) at (7.0,1.0) {\large$\gamma$};
  \node (g5) at (7.25,0.75) {\large$\gamma$};
  \node (g6) at (7.5,0.5) {\large$\gamma$};
  \node (g7) at (7.75,0.25) {\large$\gamma$};
  \node (sdirk) at (7,-0.5) {SDIRK};
\end{tikzpicture}
  \vspace{-0.2cm}
  \caption{Sparsity structure of $A$ matrices.}
  \label{fig:rk_scheme_structures}
  \vspace{-0.7cm}
\end{figure}
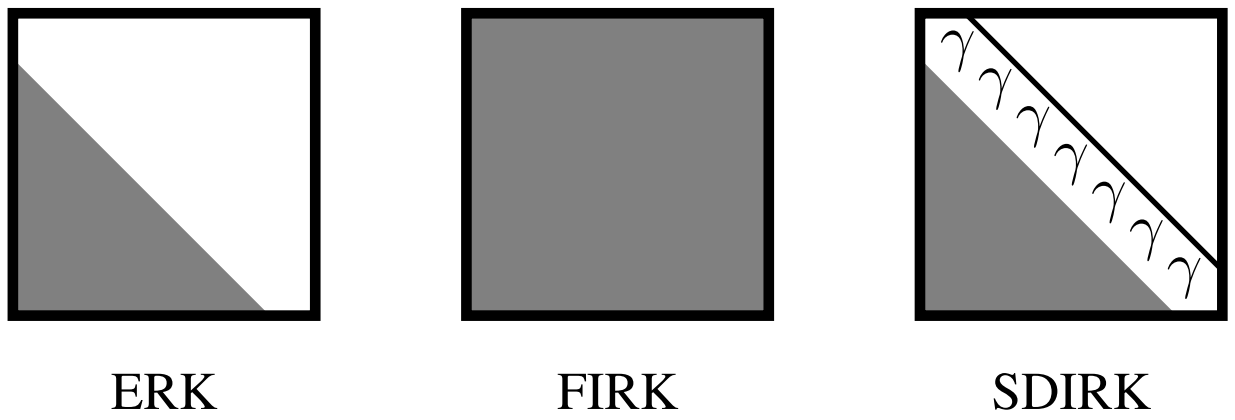

\paragraph*{Contributions}
In this paper, we discuss the efficient implementation of fixed-step-size SDIRK integrators in the open-source optimal control and MPC software {\acados}.
The implementation provides the capability for simulation as well as forward sensitivities, adjoint sensitivities, and symmetric Hessian propagation~\cite{Quirynen2017}.
We discuss inexact Newton schemes and structure-exploiting calculations which allow for significant theoretical, and, as we verify, practical, speedups over FIRK methods of equivalent order.
Finally, we compare our implementation on several system models and hardware platforms, both of which show the practical benefits of high-order SDIRK methods.
In particular, we evaluate our implementation on compute constrained embedded hardware typical of NMPC applications.

\section{Background}
In this section, we first provide an overview of Runge-Kutta methods and their implementation.
For ease of exposition, in the rest of this paper we focus only on ODEs {$\xdot = f(x,u)$}.
Our implementation additionally supports DAEs \eqref{eq:dae}.
We then review some background on dense numerical linear algebra, which we use for our theoretical analysis.
\subsection{Runge-Kutta Methods}
The classical differential Runge-Kutta form for a single fixed step size $\Delta t$ can be written as:
\begin{subequations}
  \label{eq:rk_diff}
  \begin{align}
    k_i &= f(x_n + \Delta t\sum_{j=1}^{\ns}a_{ij}k_{j}, u_n),\quad i=1,\ldots,\ns,\\
    x_{n+1} &= x_n + \Delta t\sum_{i=1}^{\ns}b_ik_i,
  \end{align}
\end{subequations}
with a Butcher tableau written as:
\begin{equation*}
  \begin{array}{c|c}
    c&A\\\hline
    & b
  \end{array}
  =
  \begin{array}{c|cccc}
    c_1&a_{11}&a_{12}&\cdots&a_{1\ns}\\
    c_2&a_{21}&a_{22}&\cdots&a_{2\ns}\\
    \vdots&\vdots&\vdots&\ddots&\vdots\\
    c_{\ns}&a_{\ns 1}&a_{\ns 2}&\cdots&a_{\ns \ns}\\\hline
    &b_1&b_2&\cdots&b_{\ns}
  \end{array}
\end{equation*}
where $A\in\R^{\ns\times\ns}$, $b,c\in \R^{\ns}$, and $k_i\in\R^\nstate$ for $i=1,\ldots,\ns$ are differential integration variables.
When the $\ns$ abscissae $c_i$ are all unique (that is, the method is \textit{non-confluent}~\cite{Kennedy2016}), the equations in \eqref{eq:rk_diff} can be rewritten in the so-called \textit{integral} form:
\begin{align}
  \label{eq:rk_integral}
  \Delta tf(s_i,u_n) &= \sum_{j=1}^{\ns}d_{ij}(s_{j} - x_n),\quad i=1,\ldots,\ns,\\
  x_{n+1} &= x_n + \sum_{i=1}^{\ns}\tilde{b}_i(s_i-x_n),
\end{align}
with
\begin{equation*}
  \begin{array}{c|c}
    c&D\\\hline
    \rule{0pt}{1.0\normalbaselineskip} & \tilde{b}
  \end{array}
\end{equation*}
where $D\in\R^{\ns\times\ns} = \inv{A}$ and $\tilde{b}\in \R^{\ns} = \transp{D}b$~\cite{Rawlings2017}.
In this form $s_i\in\R^\nstate$ for $i=1,\ldots,\ns$, are the values of the differential state at the abscissa $c_i$, and $s\in\R^{\ns\nstate}$ is the concatenation of these values.

In the rest of this paper, we discuss DIRK methods, that is, methods where the matrix $A$ (and therefore also $D$) is lower triangular.
Several extensions of the DIRK concept exist, of which we focus on the Singly Diagonally Implicit Runge-Kutta (SDIRK) methods.
SDIRK methods are defined by a lower triangular tableau, with equal diagonal elements $a_{ii} = \gamma$, with $i=1,\ldots,\ns$, as shown in \Cref{fig:rk_scheme_structures}.
This property additionally holds for the integral form with $d_{ii} = \inv{\gamma}$, $i=1,\ldots,\ns$.

In both the differential \eqref{eq:rk_diff} and integral \eqref{eq:rk_integral} case, the RK integration step reduces to the solution of a square system of nonlinear equations of dimension $\ns\nstate$.
In the integral case this is written:
\begin{equation}
  \label{eq:rk_g}
  h_i(s;x_n,u_n) = \sum_{j=1}^{\ns}d_{ij}(s_{j} - x_n) - hf(s_i,u_n) = 0,
\end{equation}
for $i=1,\ldots,\ns$.
In a practical implementation, this involves applying Newton's method to the system
\begin{subequations}
  \label{eq:newton_system}
  \begin{align}
    H^k\Delta s &= -h^k\\
    s^{k+1} &= s^k + \Delta s,
  \end{align}
\end{subequations}
where $H^k = \nabla_sh(s^k;x_n,u_n)\in\R^{\ns\nstate\times\ns\nstate}$, and $h^k = h(s^k;x_n,u_n)\in\R^{\ns\nstate}$.

\subsection{Numerical Linear Algebra}
There are two major computational costs associated with implicit RK schemes: function (and Jacobian) evaluation, and the numerical linear algebra involved in the (possibly inexact) Newton's method \eqref{eq:newton_system} used to solve the corresponding set of nonlinear equations \eqref{eq:rk_g}.
We measure these costs as a number of Floating Point Operations (FLOPs), e.g. additions, multiplications, and divisions.
In general a linear system
\[
  Ax = b,
\]
is solved for $x\in\R^{n}$ by first factorizing the matrix {$A\in\R^{n\times n}$}, most commonly by the LU factorization $A = LU$ with $L$ lower triangular, and $U$ upper triangular and (in the most naive case) done in $\frac{2}{3}n^3$ FLOPs.
These factors are then used in sequence to solve:
\begin{equation*}
  Ly=b,\quad Ux=y,
\end{equation*}
with $y\in\R^{n}$.
Each of these operations, called triangular back-solves, take a total of $2n^2-n$ FLOPs.
Note that if we are solving the matrix equation $AX = B$ with ${X,B\in\R^{n\times m}}$, the corresponding back-solves require $2mn^2-mn$ FLOPs.
A general dense matrix multiplication of $AB$ with {$A\in\R^{n\times p}$}, {$B\in\R^{p\times m}$}, takes $nm(2p-1)$ FLOPs.
Finally, when $A\in\R^{n\times n}$ is diagonal the multiplication with $B\in \R^{n\times m}$ only requires $n(2m-1)$ FLOPs~\cite{Duff1987}.

\section{Exploiting Structure}
In this section, we discuss the structure of the Newton system \eqref{eq:newton_system}, which is induced by the assumptions in DIRK methods.
To this end, we first discuss the computational costs involved in forward simulation of the system, followed by those involved in the calculation of sensitivities.
We discuss how the block diagonal structure of DIRK methods can be exploited, as well as exploiting the common elements on the diagonal of $A$ in SDIRK methods.
We then evaluate the theoretical improvement of the different formulations in terms of FLOP counts, when compared to FIRK methods.

\subsection{Computational Complexity for FIRK}
While FLOPs are not a complete evaluation of a numerical algorithm, there is often strong correlation between wall time performance and the amount of work an algorithm does.
Therefore, we first evaluate the computational cost for a fixed step size FIRK method with $\ns$ stages.
The computation required can be split into two components: the solution of $h(s,x,z,u)=0$ via Newton's method, and sensitivity computation via application of the implicit function theorem.

During Newton's method there are two operations which constitute the vast majority of the computational effort: factorizing the matrix $H$, and back-solving for $\Delta s$.
In the case of an $\ns$ stage FIRK method this corresponds to:
\begin{align*}
  \label{eq:c_newton_firk}
  \cnewton &= \cfactor + \cbacksolve\\
                     &= \frac{2}{3}(\ns\nstate)^3 + 2(\ns\nstate)^2 - 2(\ns\nstate),
\end{align*}
which must be done at each Newton iteration.
The sensitivity calculation phase then consists of a matrix valued back-solve:
\begin{equation}
  \label{eq:sens_calc}
  S_{xu} = \inv{H}
  \begin{bmatrix}
    \nabla_{x_n}h\\
    \nabla_{u}h
  \end{bmatrix},
\end{equation}
where $S_{xu}$ is the concatenation of matrix $A$ and $B$ in \eqref{eq:lin_dyn}.
If the factors of the matrix $H$ are already available at the solution, then this corresponds to a matrix valued back-solve with the FLOP cost:
\begin{equation*}
  \csens = (\ncontrol+\nstate)(2(\ns\nstate)^2 - 2(\ns\nstate)).
\end{equation*}
The total cost in FLOPs for a single step of a FIRK integrator is on the order of:
\begin{align*}
  \ctotal &= M\cnewton + \csens\\
  &\approx \frac{2}{3}M\nstate^3\ns^3 + 2(\nstate^3 + (\ncontrol+1)\nstate^2)\ns^2,
\end{align*}
with $M$ being the number of Newton iterations required for convergence to a prescribed tolerance.
In practice, we observe that performing a few Newton iterations suffices to achieve a low tolerance, practically speaking 3 to 5.

\subsection{Lower Block Triangular Form}
When the matrix $A$ (or $D$) is lower triangular, the Jacobians of both the differential and integral Runge-Kutta equations are block lower triangular:
\begin{equation*}
  H^k =
  \begin{bmatrix}
    H^k_{11}& 0 & \cdots & 0\\
    H^k_{21}& H^k_{22} & \cdots & 0\\
    \vdots& \vdots & \ddots & \vdots\\
    H^k_{\ns 1}& H^k_{\ns 1} & \cdots & H^k_{\ns \ns}
  \end{bmatrix},
\end{equation*}
with each block $H^k_{ij}\in\R^{\nstate\times\nstate}$.
It is well known~\cite{Duff1987} that block triangular systems of the form in \Cref{eq:newton_system} can be solved via the recursion:
\begin{equation}
  \label{eq:backsolve}
  H^k_{ii}\Delta s_{i} = -s^k_{i} - \sum_{j=1}^{i-1}H^k_{ij}\Delta s_j, \; \text{for} \; i=1,\ldots,\ns.
\end{equation}

When using this structured solve the FLOP cost of the factorization is reduced to
\begin{equation*}
  \cfactorhat = \frac{2}{3}\ns \nstate^3,
\end{equation*}
as it is only necessary to factor the $\ns$ square matrices $H_{ii}$ of size $\nstate\times\nstate$. 
In the case of the formulation \eqref{eq:rk_integral}, we further observe that the block matrices of $H^k$ below the diagonal are all scaled identity matrices.
This simplifies \eqref{eq:backsolve} to
\begin{equation*}
  H^k_{ii}\Delta s_{i} = -s^k_{i} - \sum_{j=1}^{i-1}a_{ij}\Delta s_j,
\end{equation*}
meaning the total work needed for the back-solve reduces to
\begin{equation*}
  \cbacksolvehat = \ns(2\nstate^2 -\nstate) + \frac{1}{2}(\ns-1)\ns\nstate.
\end{equation*}
Finally, the matrix-valued back-solve required for sensitivity calculation requires $\nstate+\ncontrol$ vector back-solves as before.
In total, this yields a total amount of FLOPS of
\begin{align*}
  \ctotalhat = &M(\tfrac{2}{3}\ns n_x^3) +\\
               &(M+\nstate+\ncontrol)(\ns(2\nstate^2 -\nstate) + \tfrac{1}{2}(\ns-1)\ns\nstate).
\end{align*}

Note that if one implements the integrator using the differential form as in \eqref{eq:rk_diff}, the blocks $H^k_{ij}$ with $i\neq j$, become dense.
This means that the back-solve requires $O(\ns^2)$ dense matrix-vector multiplications of size $\nstate$, making the high-order term increase from $O(\ns\nstate^3)$ to $O(\ns^2\nstate^3)$.
This significantly impacts the practical performance of DIRK methods, making them nearly entirely outclassed for moderate amounts of stages when adjusted for equivalent order accuracy.

\subsection{Jacobian Block Reuse}
In the case of SDIRK methods, additional FLOP count reductions can be achieved in the regime where the Jacobian of the dynamics changes slowly~\cite{Christensen2025}.
When this is the case, the singly implicit characteristic can be exploited by replacing the matrices $H^k_{ii}$ with an approximation
\begin{equation*}
  \hat{H}^k = \gamma I - \Delta t \nabla_x f(x_n, u).
\end{equation*}
The inexact Newton system solve for $\Delta s$ then becomes:
\begin{equation*}
  \hat{H}^k\Delta s_{i} = -s^k_{i} - \sum_{j=1}^{i-1}a_{ij}\Delta s_j.
\end{equation*}

This reduces the number of $\nstate\times\nstate$ matrices which need to be factorized at each Newton iteration from $\ns$ to~$1$.
In practice, this reduction is significant as high-order accuracy DIRK methods require comparatively many stages.
This approximation also reduces the number of times the Jacobian of the dynamics needs to be evaluated by the same factor, which is, however, not a dominating factor in the domain of interest.

\subsection{Sensitivity Calculation Considerations}
This Jacobian reuse in the prior section is, however, only practical for reducing the computational cost of the Newton iterations.
In optimal control, we additionally require accurate forward sensitivities.
For the sensitivity calculation $S_{xu}$, we calculate the true matrix $H$ at the solution, and use the accurate factorization when applying the implicit function theorem to~\eqref{eq:rk_integral}.
This constitutes additional factorizations of $\ns$ square $\nstate\times\nstate$ matrices.
In principle if only approximate sensitivities are required, this additional cost can also be avoided, but this exacerbates inaccuracies at large step sizes.
The effects of this approximation in optimal control has been recently studied by Christensen and J{\o}rgensen~\cite{Christensen2024}, and were found to have a negative impact on achievable accuracy.

While Hessian approximations are common in practical methods for NMPC, exact Hessian methods might result in better convergence properties~\cite{Nocedal2006} and are required for computing solution sensitivities~\cite{Frey2026a}.
As a consequence, a fully featured integrator needs to provide exact Hessian propagation, requiring adjoint sensitivity calculation which corresponds to a vector back-solve with $\transp{H}$~\cite{Quirynen2017}.
While we do not discuss them in detail, the proposed implementation supports efficient computation of adjoint and second-order sensitivities by exploiting the block triangular $H$ matrix.

\begin{figure}[t!]
  \centering
  \includegraphics[width=0.85\linewidth]{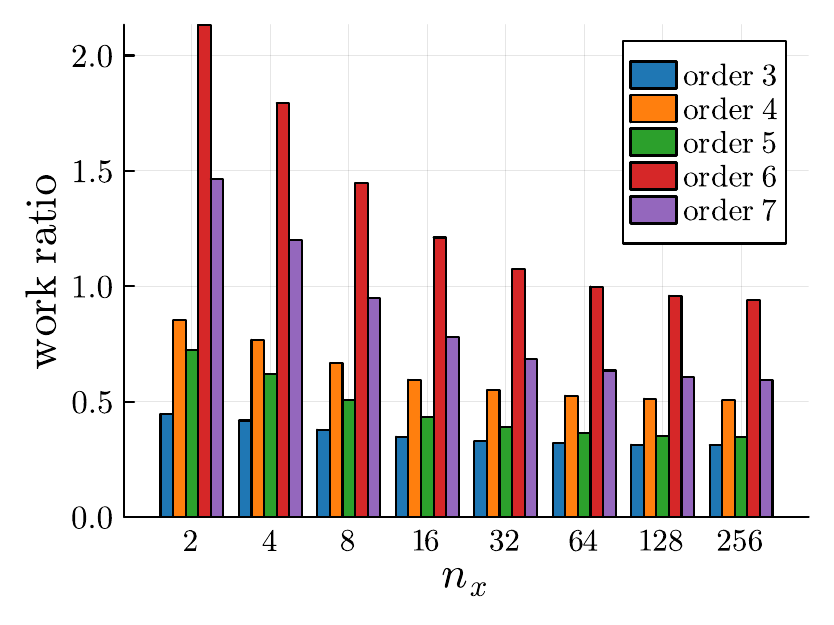}
  \vspace{-0.2cm}
  \caption{Relative amount of work, including exact forward sensitivities, for given order SDIRK versus FIRK method.}
  \label{fig:work_ratios}
  \vspace{-0.75cm}
\end{figure}

We plot approximate FLOP counts in \Cref{fig:work_ratios}, comparing SDIRK and FIRK methods of equivalent order accuracy, assuming $M=5$.
At moderate system sizes, we see significant, but not order-of-magnitude, reductions in the amount of work that needs to be done.
For DIRK methods where we cannot reuse the stage factorization, the reduction is not as significant and requires larger systems for significant benefits to be seen.
As such, we focus primarily on evaluating the available SDIRK methods.
\section{Numerical Experiments}
\begin{figure*}[t!]
  \centering
  \vspace{1mm}
  \begin{subfigure}{0.39\linewidth}
    \centering
    \includegraphics[width=\textwidth]{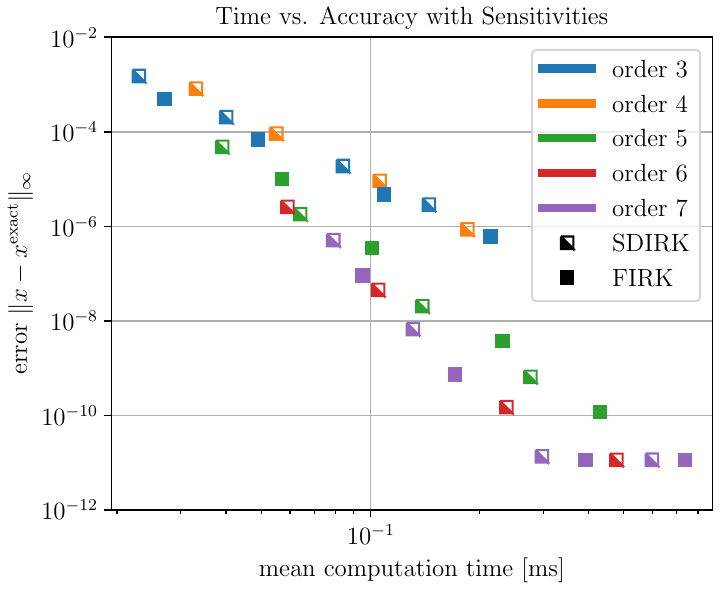}
  \end{subfigure}
  \begin{subfigure}{0.39\linewidth}
    \centering
    \includegraphics[width=\textwidth]{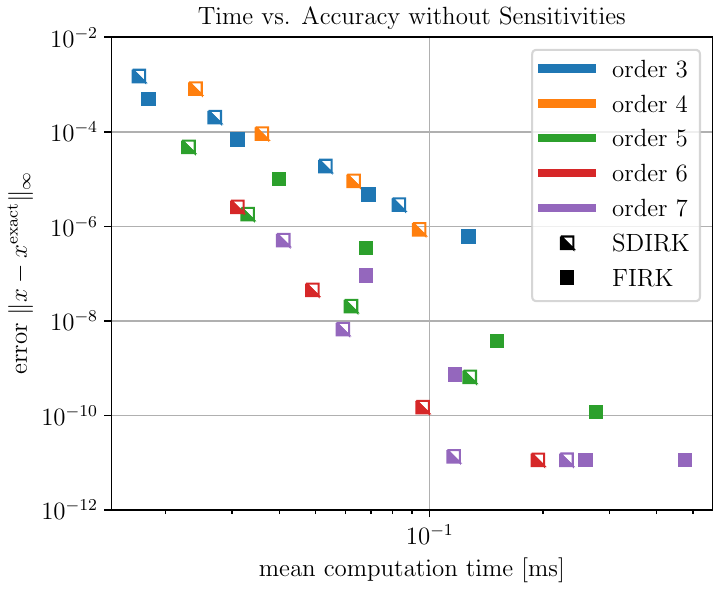}
  \end{subfigure}
  \vspace{-0.15cm}
  \caption{Pareto plots for the fixed wing aircraft model on \texttt{x86-64} CPU.}
  \label{fig:fwa_amd64}
  \vspace{-0.7cm}
\end{figure*}
This section details several numerical experiments to investigate the practical usefulness of SDIRK methods\footnote{Code: \url{https://github.com/apozharski/ecc2027_sdirk_experiments}}.
First, we compare integrator efficiency, both including and excluding forward sensitivity calculation, using SDIRK and FIRK methods on a nonlinear fixed wing aircraft with $\nstate = 22$~\cite{Gros2012e}.
To demonstrate the practical benefits of SDIRK methods in the case of embedded processors, we benchmark integrator performance on a series of platforms with tighter compute constraints.
This reveals the stronger sensitivity of embedded MPC to the theoretical reduction in FLOP count compared to desktop \texttt{x86-64} processors with wider Single Instruction, Multiple Data (SIMD) pipelines.
Finally, we evaluate the effects of scaling the state space size on SDIRK performance, using the nonlinear chain of masses example.

\begin{figure}[b!]
  \centering
  \vspace{-0.5cm}
  \includegraphics[width=0.8\linewidth]{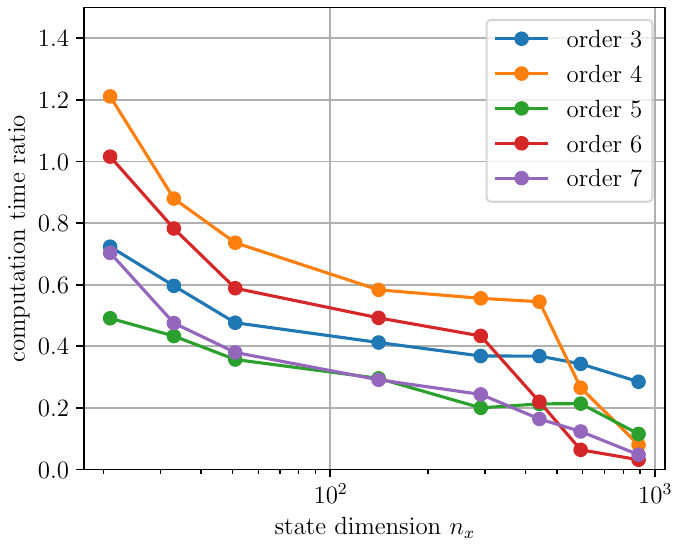}
  \vspace{-0.2cm}
  \caption{Ratio of SDIRK to FIRK computational time.
  }
  \label{fig:time_ratios}
  \vspace{-0.1cm}
\end{figure}
\subsection{Hardware and Software Configuration}
The desktop experiments are run on an Intel Core Ultra 7 155H.
In order to ensure accurate data we pin the benchmarking process to a single core and fix that core to its base frequency of 3.8 GHz.
On the desktop, we compile {\blasfeo} including the specialized kernels taking advantage of the \texttt{avx2} instruction set.

The embedded system experiments are run on an Broadcom BCM2711, quad-core Cortex-A72, clocked at 1.8 GHz, with the benchmarking process again pinned to one core.
For the embedded processor, we use the \texttt{ARMV8A\_ARM\_CORTEX\_A57} target for {\blasfeo}, which includes specialized kernels and sizing for Cortex-A72 processors, as well as taking advantage of the NEON instruction set extension.
Finally, to emulate an even more resource constrained environment we include a configuration using the Cortex-A72 as a hardware platform, paired with {\blasfeo} built against the \texttt{GENERIC} target.
This is a pure C implementation of {\blasfeo} which does not take advantage of SIMD specialized kernels, and only the number of floating point registers available~\cite{Frison2018}.

The SDIRK implementation in {\acados} that accompanies this paper implements the following SDIRK methods:
\begin{itemize}
\item An order 3 method with 2 stages\cite{Norsett1974},
\item An order 4 method with 3 stages\cite{Norsett1974},
\item An order 5 method with 5 stages\cite{Kennedy2016},
\item An order 6 stiffly accurate method with 9 stages~\cite{Alamri2024},
\item An order 7 method with 11 stages~\cite{Alamri2024},
\end{itemize}
along with several high order SDIRK methods~\cite{Alamri2024} which are uncompetitive for systems smaller than $\nstate \approx 100$.

\begin{figure*}[t!]
  \centering
  \vspace{1mm}
  \begin{subfigure}{0.39\linewidth}
    \centering
    \includegraphics[width=\textwidth]{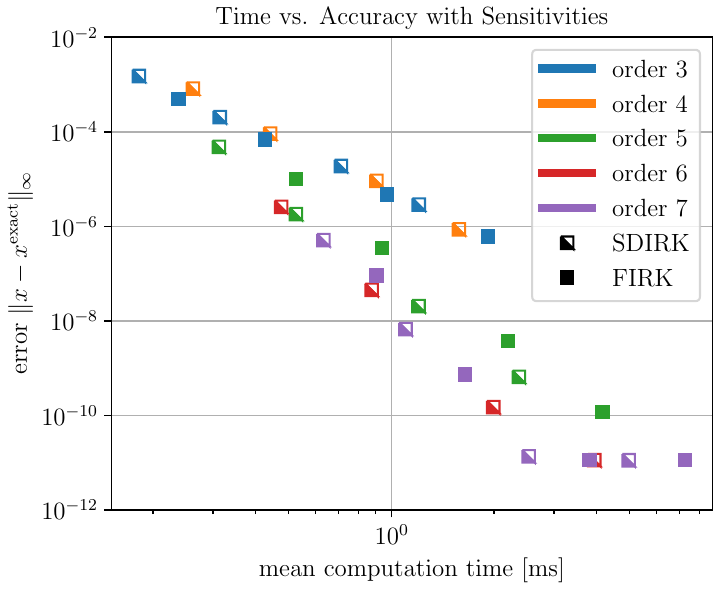}
  \end{subfigure}
  \begin{subfigure}{0.39\linewidth}
    \centering
    \includegraphics[width=\textwidth]{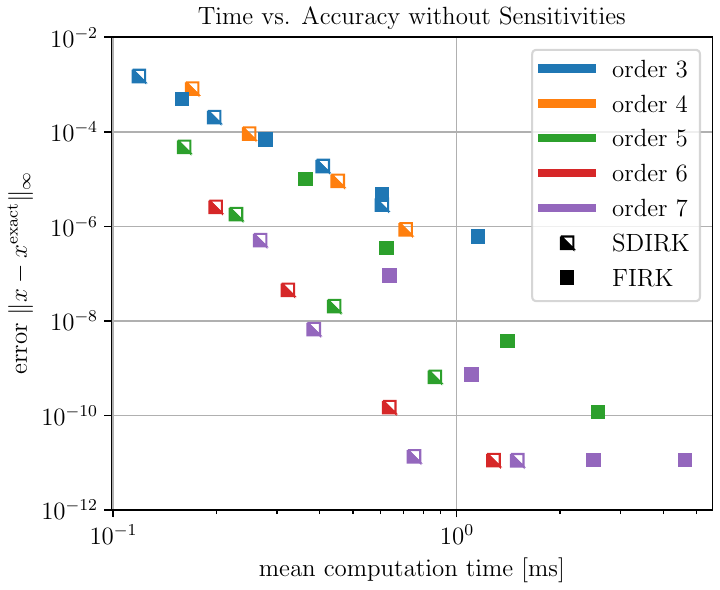}
  \end{subfigure}
  \caption{Pareto plots for the fixed wing aircraft model on \texttt{arm64} CPU.}
  \label{fig:fwa_arm64}
  \vspace{-0.3cm}
\end{figure*}
\begin{figure*}[t]
  \centering
  \begin{subfigure}{0.39\linewidth}
    \centering
    \includegraphics[width=\textwidth]{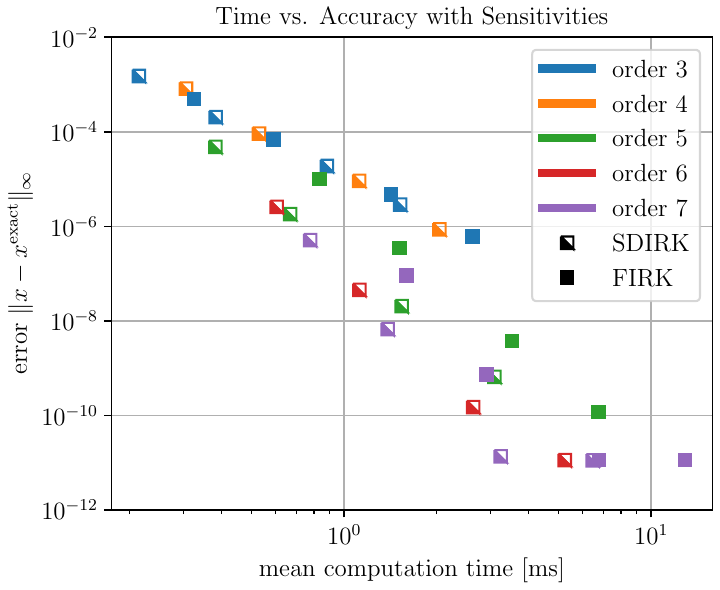}
  \end{subfigure}
  \begin{subfigure}{0.39\linewidth}
    \centering
    \includegraphics[width=\textwidth]{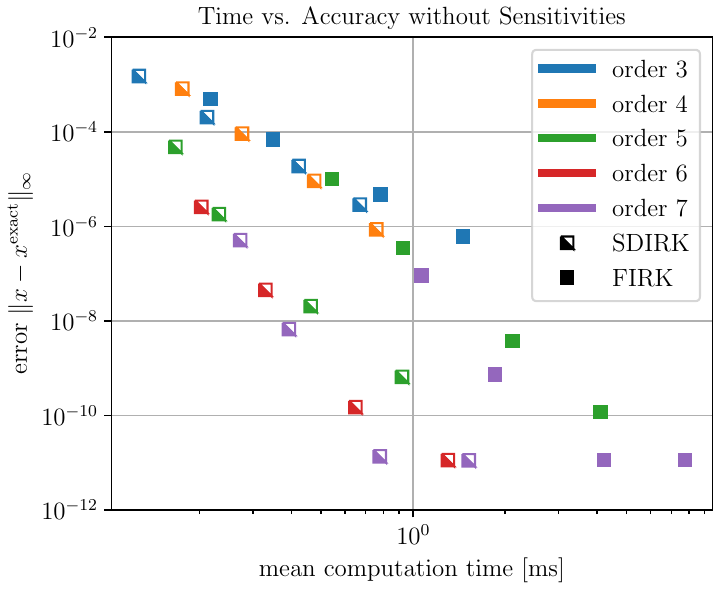}
  \end{subfigure}
  \caption{Pareto plots for the fixed wing aircraft model on \texttt{arm64} CPU, using the \texttt{GENERIC} target in {\blasfeo}.}
  \label{fig:fwa_arm64_generic}
  \vspace{-0.6cm}
\end{figure*}
\subsection{Fixed-wing Aircraft}
As a first benchmark example, we study a fixed wing aircraft with actuated control surfaces as used in~\cite{Gros2012e}, which we regard in explicit form using a symbolic inversion of the inertia matrix.
To evaluate the practical computational efficiency, we solve an initial value problem with a sequence of nominal controls.
The terminal state of this problem is then compared to a high-accuracy solution $x^{\textrm{exact}}$ computed by an 8 stage Radau-IIA method with 100 integration steps per control interval.
The control step is taken to be $200$ milliseconds, and we evaluate each integrator with $n_{\textrm{steps}} = 1,\ 2,\ 5,\ 10$ steps on each control interval.

The plots in \Cref{fig:fwa_amd64}, \Cref{fig:fwa_arm64}, and \Cref{fig:fwa_arm64_generic} show the trade-off between integrator accuracy and computation time in milliseconds, both with and without sensitivity calculations.
As we are primarily interested in optimal control applications, we focus first on the results including exact sensitivity calculation.
The resulting Pareto front on a modern laptop \texttt{x86-64} processor (\Cref{fig:fwa_amd64}) shows a relatively comparable efficiency between FIRK and SDIRK methods on a system of this size.
Evaluating the performance on the \texttt{arm64} processor with SIMD pipelines (\Cref{fig:fwa_arm64}) shows a stronger advantage to using SDIRK methods when medium to high accuracy is required.
Constraining computational resources further by not using SIMD exploiting kernels on the \texttt{arm64} processor yields the most compelling results.
In this environment, the SDIRK methods dominate the Pareto front across the entire range of accuracies and compute times.

We observe a similar improvement in relative performance in the case where sensitivities are not needed.
However, in this case, SDIRK methods dominate the Pareto front on all platforms.
Of particular interest are the very high order accuracy methods of Alamri et al.~\cite{Alamri2024} which have only recently been developed and show particularly strong efficiency at high accuracy requirements.

The above results suggest that SDIRK methods can indeed be competitive with FIRK methods in the context of NMPC.
However, there are some caveats that practitioners should be aware of when choosing methods.
Particularly, the SIMD performance of modern desktop and laptop processors make the efficiency trade-off significantly less appealing.
Additionally, we note that while much is often spoken about the order accuracy of different RK methods, little is said about truncation error.
The performance particularly of the~3 stage order 4 SDIRK method suffers significantly from its relatively large truncation error factor~\cite{Kennedy2016}.
This makes it a generally uncompetitive method even though its order accuracy is high.

\subsection{Chain of Masses}
As a second system, we regard the so-called \textit{nonlinear chain of masses} described in~\cite{Wirsching2006}.
In order to evaluate the performance of the SDIRK methods as the system size grows, we consider the system for different number of masses ranging from $3$ to $125$.

We compare the resulting runtime ratio between SDIRK and FIRK methods.
We see from these results that as the state space dimension grows, the benefits of using SDIRK methods becomes more prominent.
In practice, this is due to both the reduction in cost of the Newton steps involved in SDIRK methods, as well as cache effects as the matrix $H^k$ begins to exceed the size of the L1 cache of the CPU.
The second is quite pronounced as the state dimension scales towards $\nstate \approx 500$, where FIRK methods begin to be dominated by memory access times.
Additionally, the effects of improved throughput as matrix size increases are well documented~\cite{Frison2018}, which is shown by the worse than expected performance ratios for small systems.

\vspace{1mm}
\section{Conclusions and Future Work}
In this paper, we describe the efficient implementation of Singly Diagonally Implicit Runge-Kutta (SDIRK) methods in the {\acados} software package.
We show that when carefully implemented, these methods are on the Pareto front for efficiency with respect to the trade-off between accuracy and computational time.
This holds even when the calculation of accurate sensitivities is required.
We further show that this efficiency persists in the case of recently developed, high-order, methods with many stages~\cite{Alamri2024}.


\bibliographystyle{ieeetran}
\bibliography{syscop}
\end{document}